\documentstyle[12pt,amsfonts,amssymb,leqno]{article}

\newtheorem{proposition}{Proposition}[section]
\newtheorem{thm}[proposition]{Theorem}

\newtheorem{defn}[proposition]{Definition}

\begin{document}

\begin{center}
{\Large \bf A simple proof of
$\Sigma^1_3$ correctness of $K$} 
\end{center}
\begin{center}
\renewcommand{\thefootnote}{\fnsymbol{footnote}}
\renewcommand{\thefootnote}{arabic{footnote}}
\renewcommand{\thefootnote}{\fnsymbol{footnote}}
{\large Ralf Schindler}
\renewcommand{\thefootnote}{arabic{footnote}}
\end{center}
\begin{center} 
{\footnotesize
{\it Institut f\"ur Formale Logik, Universit\"at Wien, 1090 Wien, Austria}} 
\end{center}

\begin{center}
{\tt rds@logic.univie.ac.at}

{\tt http://www.logic.univie.ac.at/${}^\sim$rds/}\\
\end{center}

The purpose of the present paper is to present a simple proof of the following result,
which is due to John Steel.

\begin{thm}\label{thm} {\em \bf (Steel 1993, \cite[Theorem 7.9]{CMIP})} 
Let $A \subset {\mathbb R}$ be $\Pi^1_2$.
Suppose that there is some $x \in A$ such
that $(x^\dagger)^\#$ exists. Suppose also
that there is no inner model with a Woodin
cardinal. There is then an iterable lightface premouse ${\cal M}$ such
that $A \cap {\cal M} \not= \emptyset$. 
\end{thm}

We take our statement of Theorem \ref{thm} and Steel's statement of 
\cite[Theorem 7.9]{CMIP} to be basically just linguistic variants of each other.

Our proof of Theorem \ref{thm} is purely combinatorial in contrast to the proof
given in \cite[\S 7.D]{CMIP}. The latter one uses methods from descriptive set theory,
for instance the Martin-Solovay tree and the Kunen-Martin theorem. We shall be able to
avoid any serious use of descriptive set theory, except for Shoenfield absoluteness.
We believe that the argument to follow might help showing the right 
correctness results
for higher core models. 

The arguments given below would in fact enable us to prove the stronger version of 
Theorem \ref{thm} in which ``$(x^\dagger)^\#$ exists'' is replaced by 
``$x^\dagger$ exists'' (or even by something slightly less). The key problem that
still remains open, however, is
how to prove the
version of 
Theorem \ref{thm} in which ``$(x^\dagger)^\#$ exists'' is replaced
by ``$x^\#$ exists.''

As for prerequisites, an acquaintance with \cite[\S \S 1-6 and p.~58]{CMIP} will
certainly suffice. We shall also use the result of \cite{IOK}; this result would not
be needed, though, and could be replaced by a use of the weaker result \cite[Lemma
7.13]{CMIP} at the cost of introducing just a bit more notational fog to the argument
to follow. 

We'll need only a few definitions before we can commence with proving Theorem
\ref{thm}.  

\begin{defn}\label{defn-daggers}
Let $x$ be a real. An $x$-premouse ${\cal D}$ is called an {\em $x$-dagger} provided
that
for all $\xi \leq {\cal D} \cap {\rm OR}$ do we have that
$\xi = {\cal D} \cap {\rm OR}$ if and 
only if $E_\xi^{\cal D} \not= \emptyset$ and
there is some $\nu < \xi$ with $E_\nu^{\cal D} \not= \emptyset$ and
$\nu \geq {\rm crit}(E_\nu^{\cal D})^{+{\cal D}}$.
${\cal D}$ is called a {\em dagger} if ${\cal D}$ is an $x$-dagger for some real $x$. 
\end{defn}

%Here, ${\cal D}||\xi$ denotes ${\cal D}$ being cut off at $\xi$. 
${\cal D}$ is thus a
dagger if and only if ${\cal D}$ is a premouse built over a real and
${\cal D}$ is the least initial segment of itself which has two
active extenders. Note, however, that we do not require a dagger to be iterable.
Therefore, $x^\dagger$ is a dagger for any real $x$, but not the other way round.

\begin{defn}\label{defn-daggers2}
Let ${\cal D}$ be a dagger. 
We shall denote by $\kappa^{\cal D}$, $\Omega^{\cal D}$ the critical points
of the two active extenders 
of ${\cal D}$, where we understand that $\kappa^{\cal D} < \Omega^{\cal D}$. 
\end{defn}

\begin{defn}\label{defn-daggers3}
Let ${\cal D}$ be a dagger. Set $\Omega =
\Omega^{\cal D}$. Suppose that ${\cal D} \models$ ``there is no transitive model of
${\sf ZFC}$ and of height $\Omega$ which contains a Woodin cardinal.'' We then let
$K^{\cal D}$ denote Steel's core model of height $\Omega$, as being constructed inside
${\cal D}$.
\end{defn}

\cite[\S \S 1-5]{CMIP} give the recipe for how to construct $K^{\cal D}$. We remark
that the ultrapower of ${\cal D}$ by its top extender doesn't have to be well-founded
for $K^{\cal D}$ to provably exist.

\begin{defn}\label{defn-daggers4}
Let ${\cal A}$ be a transitive model of ${\sf ZFC}$. Then by $K^{\cal A}$ we denote
the model which is recursively constructed inside ${\cal A}$ in the manner of
\cite[\S 6]{CMIP}, if it exists (otherwise we let $K^{\cal A}$ undefined). 
\end{defn}

If ${\cal D}$ is a 
premouse then we let ${\cal D}||\xi$ denote ${\cal D}$ being cut off at $\xi$.
If ${\cal D}$ is as in Definition
\ref{defn-daggers3} then $K^{\cal D}$ in the sense of Definition
\ref{defn-daggers3} is identical with $K^{{\cal D}||\Omega^{\cal D}}$ in the sense of 
Definition \ref{defn-daggers4}. This follows from \cite[\S 6]{CMIP}. The two notations
introduced by definitions \ref{defn-daggers3} and \ref{defn-daggers4} cannot be
confused, as no dagger is a model of the power set axiom.

We now turn to our proof of Theorem \ref{thm}. 

\bigskip
{\sc Proof} of Theorem \ref{thm}. Fix $A$ and $x$. 
Let $A = \{ y \in {\mathbb R} \colon \Phi(y) \}$, where $\Phi(-)$ is
$\Pi^1_2$. By the hypotheses, we know that $K^{L[x^\dagger]}$
exists (cf.~\cite[p.~58]{CMIP}). 
Let us write $K = K^{L[x^\dagger]}$. 

There is a tree $T \in K$ of height
$\omega$ searching for a quadruple $(y,{\cal D},{\cal T},\sigma)$ 
with the properties that:

\bigskip
\noindent $\bullet \ $ $y$ is a real,

\noindent $\bullet \ $ ${\cal D}$ is a $y$-dagger with ${\cal D} \models \Phi(y)$, 

\noindent $\bullet \ $ ${\cal T}$ is an iteration tree on $K$ of countable
successor length, and

\noindent $\bullet \ $ $\sigma \colon K^{\cal D} \rightarrow {\cal M}_\infty^{\cal T}
|| \alpha$ is elementary for some $\alpha \leq {\cal M}_\infty^{\cal T} \cap {\rm OR}$.

\bigskip
We leave it to the reader's discretion to construct such a tree $T$.

Let us write $p[T] = \{ y \colon \exists {\cal D} \ \exists {\cal T} \ \exists \sigma \ 
(y,{\cal D},{\cal T},\sigma) \in [T] \}$. We claim that $\emptyset \not= p[T] \subset
A$. This will establish Theorem \ref{thm}. 

\bigskip
{\bf Claim 1.} $[T] \not= \emptyset$ (in $V$, and hence in $K$).

\bigskip
{\sc Proof.} 
Set ${\cal D} = x^\dagger$. Let $({\cal U},{\cal T}')$ denote the
coiteration of $K^{\cal D}$ with $K$, which exists inside $L[x^\dagger]$. 
(We here use the fact that $K^{\cal D}$ is iterable in $L[x^\dagger]$.) 
We'll have that $\pi_{0 \infty}^{\cal U} \colon
K^{\cal D} \rightarrow {\cal M}_\infty^{{\cal T}'} || \alpha$ for some $\alpha \leq 
{\cal M}_\infty^{{\cal T}'} \cap {\rm OR}$. 
However, as ${\cal T}'$ might be uncountable, we'll have to take a Skolem hull
to finish the argument. 

Let $\tau \colon {\bar H} \rightarrow
H_\theta$, where $\theta$ is regular and large enough, ${\bar H}$ is countable and
transitive, and $\{ x , {\cal D}, {\cal U} , {\cal T}' \} \subset {\rm ran}(\tau)$.
Let us copy $\tau^{-1}({\cal T}')$ 
%(which is a tree on $\sigma^{-1}(K)$) 
onto $K$,
using $\tau$. We get a countable 
tree $\tau^{-1}({\cal T}')^{\tau}$ on $K$; let us write ${\cal T}$ for
$\tau^{-1}({\cal T}')^{\tau}$. We also get 
a last copy map
$\varphi \colon {\cal M}_\infty^{\tau^{-1}({\cal T}')} \rightarrow {\cal
M}_\infty^{\cal T}$. As $\tau \upharpoonright {\cal D} \cup \{ {\cal D} \} = {\rm
id}$, we then have that
$$\varphi \circ \tau^{-1}(\pi^{\cal U}_{0 \infty}) \colon K^{\cal D} \rightarrow
{\cal M}_\infty^{\cal T} || \varphi(\alpha) {\rm , }$$
where we understand that $\varphi(\alpha) = {\cal M}_\infty^{\cal T} \cap {\rm
OR}$ if $\alpha = {\cal M}_\infty^{{\cal T}'} \cap {\rm OR}$ (a case which actually
never comes up).
Setting $\sigma = \varphi \circ \tau^{-1}(\pi^{\cal U}_{0 \infty})$, 
we'll thus have that $(x,{\cal D},{\cal T},
\sigma) \in [T]$.

\hfill $\square$ (Claim 1)

\bigskip
{\bf Claim 2.} $p[T] \subset A$.

\bigskip
{\sc Proof.} Let $(y,{\cal D},{\cal T},
\sigma) \in [T]$. 
Let $\nu < {\cal D} \cap {\rm OR}$ be such that $E_\nu^{\cal D} \not= \emptyset$, 
${\rm crit}(E_\nu^{\cal D}) = \kappa^{\cal D}$, and
$\nu \geq (\kappa^{\cal D})^{+{\cal D}}$. Let us write $E = E_\nu^{\cal D}$.
By Shoenfield absoluteness it will suffice to prove that ${\cal D}
|| \Omega^{\cal D}$
is iterable by $U$ and its images. 
Let $({\cal D}_i , \pi_{ij} \colon i \leq j \leq \gamma)$ be a putative iteration of
${\cal D} || \Omega^{\cal D}$, where 
%${\cal D}_0 = {\cal D} || \Omega^{\cal
%D}$, 
$\gamma < \omega_1$ and 
$\pi_{i i+1} \colon {\cal D}_i
\rightarrow_{\pi_{0 i}(E)} {\cal D}_{i+1}$ for all $i < \gamma$. We have to prove that
${\cal D}_\gamma$ is well-founded.
For $i < \gamma$ let us write $E_i$ for $\pi_{0i}(E)$.

Let us first assume that $\gamma$ is a successor ordinal, $\gamma = \delta +1$, say.
Then ${\cal D}_\gamma$ is obtained by an 
internal ultrapower of ${\cal D}_\delta$. We may thus argue inside ${\cal D}_\delta$
to conclude that ${\cal D}_\gamma$ is well-founded.
 
Let us now assume that $\gamma$ is
a limit ordinal. $({\cal D}_i , \pi_{ij} \colon i \leq j \leq \gamma)$ is then the
direct limit of $({\cal D}_i , \pi_{ij} \colon i \leq j < \gamma)$.

We shall, for each $\delta < \gamma$, recursively construct an iteration tree ${\cal
T}^\delta$ of length $\beta_\delta +1$ on $K^{\cal D}$, 
and we shall inductively verify that the
following clauses hold true:

\bigskip
(a)${}_\delta \ $ ${\cal
T}^i = {\cal
T}^{\delta} \upharpoonright \beta_i +1$ for all $i \leq
\delta$,

(b)${}_\delta \ $ ${\cal M}_{\beta_\delta}^{{\cal T}^\delta} = K^{{\cal D}_\delta}$,
and

(c)${}_\delta \ $ $\pi_{\beta_i 
\beta_\delta}^{{\cal T}^\delta} = \pi_{i \delta} \upharpoonright
K^{{\cal D}_i}$ for all $i \leq \delta$.

\bigskip
However, this is a straightforward task. 
To get started, let us apply
\cite[Corollary 3.1]{IOK} inside ${\cal D}$ to get an iteration tree ${\cal U}$ on
$K^{\cal D}$ with ${\cal M}_\infty^{\cal U} = K^{{\rm Ult}({\cal D};E)}$ and
$\pi_{0 \infty}^{\cal U} = \pi_{01} \upharpoonright K^{\cal D}$. Notice that we
may expand the model ${\cal
D} || \Omega^{\cal D}$ by a predicate coding ${\cal U}$, which we shall also denote by
${\cal U}$, to get $({\cal
D} || \Omega^{\cal D};{\cal U})$ as an amenable model. We may and shall construe
$({\cal D}_i , \pi_{ij} \colon i \leq j \leq \gamma)$ as an iteration of 
$({\cal
D} || \Omega^{\cal D};{\cal U})$ rather than of ${\cal
D} || \Omega^{\cal D}$. For $i < \gamma$ 
we'll write $\pi_{0 i}({\cal U})$ for the image of ${\cal U}$
under $\pi_{0 i}$, which is well-defined by the amenability of $({\cal
D} || \Omega^{\cal D};{\cal U})$.
We'll have that $${\cal D}_0 = ({\cal
D} || \Omega^{\cal D};{\cal U}) \ \models \ {\rm `` }{\cal M}_\infty^{\cal U} = K^{{\rm
Ult}(V;E)}{\rm , \ and }$$ $$ \ \ \ \ \ \ \ \ \ \ \ \ \ \ \ \ \ \ \ \ \ \ \ 
\pi_{0 \infty}^{\cal U} = \pi_E
\upharpoonright K {\rm ."}$$

Let us now construct $({\cal T}^\delta \colon \delta < \gamma)$. 
To commence, we let ${\cal T}^0$ be trivial. (a)${}_0$, 
(b)${}_0$, and (c)${}_0$ are trivially true.
Now suppose that ${\cal T}^\delta$ has been constructed for some $\delta < \gamma$.
We may then simply let ${\cal T}^{\delta +1}$ be the concatenation of ${\cal T}^\delta$
with $\pi_{0 \delta}({\cal U})$. The elementarity of the map
$\pi_{0 \delta}$
gives that $${\cal D}_\delta \ \models \ {\rm `` }
{\cal M}_\infty^{\pi_{0 i}({\cal U})} = K^{{\rm Ult}(V;E_i)}{\rm , \ and }$$
$$\pi_{0 \infty}^{\pi_{0 i}({\cal U})} = \pi_{E_i}
\upharpoonright K {\rm ."}$$
(a)${}_{\delta +1}$, (b)${}_{\delta +1}$, and 
(c)${}_{\delta +1}$ will then be evident.
Finally, let $\delta < \gamma$ be a limit ordinal and suppose
that ${\cal
T}^i$ has been constructed for every $i < \delta$.
Let
${\cal T}'$ be the ``union'' of all ${\cal
T}^i$ for $i < \delta$, and
let $b$ be the unique cofinal branch through ${\cal T}'$ which is generated by $\{
\beta_i \colon i < \delta \}$. 
As (b)${}_i$ and (c)${}_i$ hold 
for all $i < \delta$
we'll have that $K^{{\cal D}_\delta} = {\cal M}^{{\cal T}'}_b$ and $\pi_{\beta_i
b}^{{\cal T}'} = \pi_{i \delta} \upharpoonright K^{{\cal D}_i}$
for all $i < \delta$.
We may thus let ${\cal T}^\delta$ be that extension of ${\cal T}'$ which adds the
branch $b$ as well as the final model ${\cal M}^{{\cal U}'}_b$. Then  
(a)${}_{\delta}$, (b)${}_{\delta}$, and 
(c)${}_{\delta}$ are evident.

We may now let ${\cal U}^*$ be the union of all ${\cal U}^\delta$ for $\delta <
\gamma$. 
Let $b$ be the unique cofinal branch through ${\cal U}^*$, which is given by $\{
\beta_\delta \colon \delta < \gamma \}$. 
The tree ${\cal T}$ on $K$ and the map $\sigma$ witness 
that $K^{\cal D}$ is iterable (in
$L[x^\dagger]$, and hence in $V$). The model ${\cal
M}_b^{{\cal U}^*}$ is thus well-founded.  
As (c)${}_{\delta}$ holds for all $\delta < \gamma$, 
we now have an $\in$-isomorphism between the ordinals of ${\cal D}_\gamma$ and the
ones of ${\cal
M}_b^{{\cal U}^*}$. Therefore, ${\cal D}_\gamma$ is well-founded, too.

\hfill $\square$ (Claim 2)

\bigskip
Now let $y \in p[T] \cap K$. Let $\epsilon$ least such that 
$y \in K || \epsilon +1$. Then $K || \epsilon +1$
is iterable in $L[x^\dagger]$, and hence in $V$. We have found an iterable premouse as
desired.

\hfill $\square$ (Theorem \ref{thm})

\end{document}